\documentclass[a4paper,12pt]{amsart}
\usepackage{amssymb}
\usepackage{cite}
\usepackage{ifthen}
 \usepackage[dvips]{graphicx}
\nonstopmode \numberwithin{equation}{section}
\usepackage{amssymb}
\usepackage{ifthen}
\usepackage{graphicx}
\usepackage{amsmath}
\usepackage[T1]{fontenc} 
\usepackage[utf8]{inputenc}
\usepackage[usenames,dvipsnames]{color}
\usepackage{color}
\usepackage[english]{babel}
\usepackage{fancyhdr}
\usepackage{fancybox}
\usepackage{tikz}

\theoremstyle{plain}
\newtheorem{prop}{Proposition}

\newtheorem{conj}{Conjecture}

\theoremstyle{definition}
\newtheorem{defn}{Definition}[section]

\newtheorem{prob}{Problem}
\newtheorem{rem}{Remark}[section]
\newtheorem{ques}{Question}[section]
\newtheorem{lem}{Lemma}[section]
\newtheorem{cor}{Corollary}[section]
\newtheorem{thm}{Theorem}[section]

\newcounter{minutes}
\divide\time by 60
\newcounter{hours}
\multiply\time by 60
\addtocounter{minutes}{-\time}

\newcounter {own}
\def\theown {\thesection       .\arabic{own}}

\newenvironment{pf}[1][]{%
 \vskip 3mm
 \noindent
 \ifthenelse{\equal{#1}{}}%
  {{\slshape Proof. }}%
  {{\slshape #1.} }%
 }%
{\qed\bigskip}

\newcounter{alphabet}

\def\be{\begin{equation}}
\def\ee{\end{equation}}

\newcommand{\bee}{\begin{enumerate}}
\newcommand{\eee}{\end{enumerate}}

\newcommand{\blem}{\begin{lem}}
\newcommand{\elem}{\end{lem}}
\newcommand{\bthm}{\begin{thm}}
\newcommand{\ethm}{\end{thm}}
\newcommand{\bcor}{\begin{cor}}
\newcommand{\ecor}{\end{cor}}
\newcommand{\beg}{\begin{examp}}
\newcommand{\eeg}{\end{examp}}
\newcommand{\begs}{\begin{examples}}
\newcommand{\eegs}{\end{examples}}

\newcommand{\bdefn}{\begin{defn}}
\newcommand{\edefn}{\end{defn}}

\newcommand{\bprob}{\begin{prob}}
\newcommand{\eprob}{\end{prob}}
\newcommand{\bei}{\begin{itemize}}
\newcommand{\eei}{\end{itemize}}

\newcommand{\bcon}{\begin{conj}}
\newcommand{\econ}{\end{conj}}
\newcommand{\bcons}{\begin{conjs}}
\newcommand{\econs}{\end{conjs}}
\newcommand{\bprop}{\begin{prop}}
\newcommand{\eprop}{\end{prop}}
\newcommand{\br}{\begin{rem}}
\newcommand{\er}{\end{rem}}
\newcommand{\brs}{\begin{rems}}
\newcommand{\ers}{\end{rems}}
\newcommand{\bo}{\begin{obser}}
\newcommand{\eo}{\end{obser}}
\newcommand{\bos}{\begin{obsers}}
\newcommand{\eos}{\end{obsers}}
\newcommand{\bpf}{\begin{pf}}
\newcommand{\epf}{\end{pf}}
\newcommand{\ba}{\begin{array}}
\newcommand{\ea}{\end{array}}
\newcommand{\beq}{\begin{eqnarray}}
\newcommand{\beqq}{\begin{eqnarray*}}
\newcommand{\eeq}{\end{eqnarray}}
\newcommand{\eeqq}{\end{eqnarray*}}

\begin{document}

\title{Improved Bohr inequalities for certain classes of harmonic mappings}

\author{Molla Basir Ahamed}
\address{Molla Basir Ahamed,
Department of Mathematics,
Jadavpur University,
Kolkata-700032, West Bengal, India.}
\email{mbahamed.math@jadavpuruniversity.in}

\author{Vasudevarao Allu}
\address{Vasudevarao Allu,
	School of Basic Science,
	Indian Institute of Technology Bhubaneswar,
	Bhubaneswar-752050, Odisha, India.}
\email{avrao@iitbbs.ac.in}

\subjclass[{AMS} Subject Classification:]{Primary 30C45, 30C50, 30C80}

\def\thefootnote{}
\footnotetext{ {\tiny File:~\jobname.tex,
printed: \number\year-\number\month-\number\day,
          \thehours.\ifnum\theminutes<10{0}\fi\theminutes }
} \makeatletter\def\thefootnote{\@arabic\c@footnote}\makeatother

\begin{abstract}
The Bohr radius for the class of harmonic functions of the form $ f(z)=h+\overline{g} $ in the unit disk $ \mathbb{D}:=\{z\in\mathbb{C} : |z|<1\} $, where $ h(z)=\sum_{n=0}^{\infty}a_nz^n $ and $ g(z)=\sum_{n=1}^{\infty}b_nz^n $ is to find the largest radius $ r_f $, $ 0<r_f<1 $ such that
\begin{equation*}
	\sum_{n=1}^{\infty}\left(|a_n|+|b_n|\right)|z|^n\leq d(f(0),\partial f(\mathbb{D}))
\end{equation*}
holds for $ |z|=r\leq r_f $, where $ d(f(0),\partial f(\mathbb{D})) $ is the Euclidean distance between $ f(0) $ and the boundary of $ f(\mathbb{D}) $. In this paper, we prove two-type of improved versions of the Bohr inequalities, one for certain class of harmonic and univalent functions and the other is for stable harmonic mappings. It is observed in the paper that to obtain sharp Bohr inequalities it is enough to consider any non-negative real coefficients of the quantity $ S_r/\pi $. As a consequence of the main result, we prove a corollaries showing the precise value of the sharp Bohr radius. 
\end{abstract}

\maketitle
\pagestyle{myheadings}
\markboth{Molla Basir Ahamed and Vasudevarao Allu}{Improved Bohr inequalities for certain classes of harmonic mappings}

\section{Introduction}
It is well-known that harmonic functions are famous for their use in the study of minimal surfaces and also play important roles in a variety of problems in applied mathematics (e.g., see Choquet\cite{Choquet-1945}, Dorff \cite{Dorff-2003}, Duren \cite{Duren-Harmonic-2004} or Lewy \cite{Lew-BAMS-1936}). More precisely, harmonic mappings play the natural role in parameterizing minimal surfaces in the context of differential geometry. However, planar harmonic mappings have application not only in the differential geometry but also in the various fields of engineering, physics, operations research and other intriguing aspects of applied mathematics. The theory of harmonic functions has been used to study and solve fluid flow problems (see \cite{Aleman-2012}). The theory of univalent harmonic functions having prominent geometric properties like starlikeness, convexity and close-to-convexity appears naturally while dealing with planner fluid dynamical problems. For instance, the fluid flow problem on a convex domain satisfying an interesting geometric property has been extensively studied by Aleman and Constantin  \cite{Aleman-2012}. With the help of geometric properties of harmonic mappings, Constantin and Martin \cite{constantin-2017} have obtained a complete solution of classifying all two dimensional fluid flows.\vspace{1.2mm}

\par In $ 1914 $, Bohr \cite{Bohr-1914}, while dealing with a problem connected with Dirichlet series and number theory, proved that if the power series $ f(z)=\sum_{n=0}^{\infty}a_nz^n $ converges in the unit disk $ \mathbb{D}=\{z\in\mathbb{C} : |z|<1\} $ and $ |f(z)|<1 $, then the majorant series $ M_f(r)= \sum_{n=0}^{\infty}|a_n||z|^n $ of $ f $ satisfies 
\begin{equation}\label{e-11.1}
	M_{f_0}(r)=\sum_{n=1}^{\infty}|a_n||z|^n\leq 1-|a_0|=1-|f(0)|=d(f(0),\partial f(\mathbb{D}))
\end{equation} 
for all $ z\in\mathbb{D} $ with $ |z|=r\leq 1/6 $, where $ f_0(z)=f(z)-f(0) $ and $ d $ denotes the Euclidean distance. The interesting inequality \eqref{e-11.1} is known as the Bohr inequality. The largest $ r\leq 1 $ such that the above inequality holds is referred to as the Bohr radius for the case of unit disk. Later, Wiener, Riesz and Schur, independently established that the inequality \eqref{e-11.1} holds for $ |z|\leq 1/3 $, and hence proved that $ 1/3 $ is the best possible. The radius $ 1/3 $ is famously known as the Bohr radius. Henceforth, if there exists a positive real number $ r_0 $ such that the inequality \eqref{e-11.1} holds for every element of a class $ \mathcal{F} $ for $ 0\leq r\leq r_0 $ and fails when $ r>r_0 $, then we shall say that $ r_0 $ is an sharp bound for the inequality w.r.t. the class $ \mathcal{F} $.\vspace{1.5mm}

\par Interest in this inequality was revived in $ 1995 $ when Dixon \cite{Dixon & BLMS & 1995}, used it to settle in disproving the conjecture that a non-unital Banach algebra that satisfies the von Neumann inequality must be isometrically isomorphic to a closed subalgebra of $ B(H) $ for some Hilbert space $ H $. Initially, the problem was considered by Bohr while working on the absolute convergence of Dirichlet series of the form $ \sum a_nn^{-s} $, but in the recent years it becomes an active area of research in the geometric function theory. In fact, this theorem gets much attention as it has been applied to the characterization problem of Banach algebra satisfying the von Neumann inequality \cite{Dixon & BLMS & 1995}. Subsequently, Paulson and Singh \cite{Paulsen-PAMS-2004}, and Blasco \cite{Blasco-2010} have extended the Bohr inequality in the context of Banach algebra. Instead of the classes of analytic functions of one complex variables, Bohr radius has also been investigated by many researchers for the functions in several complex variable with different settings (see e.g. \cite{Bayart-AM-2014,Galicer-TAMS-2020}). For the extensive study of Bohr radius in multidimensional and operator valued functions, we refer to \cite{Aizn-PAMS-2000,Boas-Khavinson-PAMS-1997,Hamada-IJM-2009,Defant-AM-2012,Galicer-TAMS-2020,Lin-Liu-Pon-AMS-2023,Liu-Ponnusamy-BSM-2021,Liu-Ponnusamy-PAMS-2021} and references therein. \vspace{1.5mm}

\par A harmonic mapping in the unit disk $ \mathbb{D} $ is a complex-valued function $ f=u+iv $ of $ z=x+iy $ in $ \mathbb{D} $, which satisfies the Laplace equation $ \Delta f=4f_{z\bar{z}}=0 $, where $ f_{z}=(f_{x}-if_{y})/2 $ and $ f_{\bar{z}}=(f_{x}+if_{y})/2 $ and $ u $ and $ v $ are real-valued harmonic functions in $ \mathbb{D} $. It follows that the function $ f $ admits the canonical representation $ f=h+\bar{g} $, where $ h $ and $ g $ are analytic in $ \mathbb{D} $. The Jacobian $ J_f $ of $ f=h+\overline{g} $ is given by $ J_f(z)=|h^{\prime}(z)|^2-|g^{\prime}(z)|^2 $. We say that $ f $ is sense-preserving in $ \mathbb{D} $ if $ J_f(z)>0 $ in $ \mathbb{D} $. Consequently, $ f $ is locally univalent and sense-preserving in $ \mathbb{D} $ if, and only if, $ J_f(z)>0 $ in $ \mathbb{D} $; or equivalently if $ h^{\prime}\neq 0 $ in $ \mathbb{D} $ and the dilation $ \omega_f := \omega=g^{\prime}/h^{\prime} $ has the property that $ |\omega(z)|<1 $ in $ \mathbb{D} $. From the definition it is easy to see that class of harmonic mappings contains the class of self-analytic maps on unit disk, hence there are properties of harmonic mappings which are certainly not true for analytic functions. Therefore, it is natural to establish the Bohr type  inequalities and their possible sharped form for certain class of harmonic mappings. \vspace{1.5mm}

Bohr phenomenon can be studied in view of the Euclidean distance and in this paper, we study the same for certain classes of harmonic mappings. Before we go into details, we recall here the following concepts. Let $ f $ and $ g $ be two analytic functions in the unit disc $ \mathbb{D} $. We say that $ g $ is subordinate to $ f $ if there is a function $ \varphi $, analytic in $ \mathbb{D} $, $ \varphi(\mathbb{D})\subset \mathbb{D} $ and $ \varphi(0)=0 $ so that $ g=f\circ \varphi $. In particular, when the function $ f $ is univalent, $ g $ is subordinate to $ f $ when $ g(\mathbb{D}\subset f(\mathbb{D})) $ and $ g(0)=f(0) $ (see \cite[p. 190]{Duren-1983}). Consequently, when $ g $ is subordinate to $ f $, $ |g^{\prime}(0)|\leq |f^{\prime}(0)| $. The class of all function $ g $ subordinate to a fixed function $ f $ is denoted by $ \mathcal{S}(f) $ and $ f(\mathbb{D})=\Omega $.
\begin{defn}\cite{Abu-CVEE-2010}
	We say that $ \mathcal{S}(f) $ has Bohr phenomenon if for any $ g=\sum_{n=0}^{\infty}b_nz^n\in\mathcal{S}(f) $ and $ f=\sum_{n=0}^{\infty}a_nz^n $ there is a $ \rho^*_0 $, $ 0<\rho^*_0\leq 1 $ so that $ \sum_{n=1}^{\infty}|b_nz^n|\leq d(f(0), \partial \Omega)  $,  for $ |z|<\rho^*_0 $. Notice that $  d(f(0), \partial \Omega) $ denote the Euclidean distance between $ f(0) $ and the boundary of a domain $ \Omega $, $ \partial\Omega $. In particular, when $ \Omega=\mathbb{D} $, $  d(f(0), \partial \Omega)=1-|f(0)| $ and in this case $ \sum_{n=1}^{\infty}|a_nz^n|\leq d(f(0), \partial \Omega) $ reduces to $ \sum_{n=0}^{\infty}|a_nz^n|\leq 1 $.
\end{defn}
\par Equation \eqref{e-1.1} can be written as
\begin{equation}\label{e-1.2}
	d\left(\sum_{n=0}^{\infty}|a_nz^n|,|a_0|\right)=\sum_{n=1}^{\infty}|a_nz^n|\leq 1-|f(0)|=d(f(0),\partial (\mathbb{D})),
\end{equation}
where $ d $ is the Euclidean distance. More generally, a class $ \mathcal{F} $ of analytic functions $ f(z)=\sum_{n=0}^{\infty}a_nz^n $ mapping $ \mathbb{D} $ into a domain $ \Omega $ is said to satisfy a Bohr phenomenon if an inequality of type \eqref{e-1.2} holds uniformly in $ |z|\leq \rho_0 $, where $ 0<\rho_0\leq 1 $ for functions in the class $ \mathcal{F} $. Similar definition makes sense for harmonic functions (see \cite{Kay & Pon & Sha & MN & 2018}). Abu-Muhanna \cite{Abu-CVEE-2010} established the following result for subordination class $ \mathcal{S}(f) $ when $ f $ is univalent.
\begin{thm}\cite{Abu-CVEE-2010}
	If $ g=\sum_{n=0}^{\infty}b_nz^n\in\mathcal{S}(f) $ and $ f=\sum_{n=0}^{\infty}a_nz^n $ is univalent, then 
	\begin{align*}
		\sum_{n=1}^{\infty}|b_nz^n|\leq d(f(0), \partial\Omega),
	\end{align*}
	for $ |z|=\rho^*_0\leq 3-\sqrt{8}\approx0.17157 $, where $ \rho^*_0 $ is sharp for the Koebe function $ f(z)=z/(1-z)^2 $.
\end{thm}

\par  Like class of analytic self maps on unit disk, the Bohr's phenomenon for complex-valued harmonic mappings much not studied and a little is known about the sharp Bohr-type inequalities and their different versions (see, e.g.  \cite{Abu-CVEE-2010,Abu-JMAA-2014,Himadri-Vasu-P1,Nirupam-VVEE-2018,Nirupam-MonatsMath-2019}). For example, in $ 2016 $, Ali \emph{et al.} \cite{Ali & Abdul & Ng & CVEE & 2016} studied Bohr radius for the stralike log-harmonic mappings and obtain several interesting results.  Bohr-type inequalities for the class of harmonic $ p$-symmetric mappings and also for harmonic mappings with multiple zeros at the origin have been discussed by Huang \emph{et al.} \cite{Huang-Liu-Ponnusamy-MJM-2021}.  The Bohr radius for various classes of functions, for example, locally univalent harmonic mappings, $ K $-quasiconformal mappings, bounded harmonic functions, functions having lacunary series, have been studied extensively in \cite{Ismagilov-2020-JMAA,Kayumov-2018-JMAA,Kay & Pon & Sha & MN & 2018}. For more intriguing aspects of the Bohr phenomenon, we refer to the articles \cite{Aha-CMFT-2022,Ahamed-AMP-2021,Ahamed-CVEE-2021,Aizn-PAMS-2000,Ali-Jain-RM-2019,Alkhaleefah-PAMS-2019,Allu-JMAA-2021,Bhowmik-2018,Evdoridis-IM-2018,Kay & Pon & AASFM & 2019,Kayumov-CRACAD-2020,Liu-Ponnusamy-BMMS-2019,Ponnusamy-Vijayakumar-2022-Survey,Ponnusamy-RM-2020,Ponnusamy-HJM-2021} and references therein.\vspace{1.5mm} 

\par Let $ \mathcal{H} $ be the class of all complex-valued harmonic functions $ f=h+\bar{g} $ defined on the unit disk $ \mathbb{D} $, where $ h $ and $ g $ are analytic in $ \mathbb{D} $ with the normalization $ h(0)=h^{\prime}(0)-1=0 $ and $ g(0)=0 $. Let $ \mathcal{H}_0 $ be defined by $ 	\mathcal{H}_0=\{f=h+\bar{g}\in\mathcal{H} : g^{\prime}(0)=0\}. $ Then each $ f=h+\bar{g}\in \mathcal{H}_0 $ has the following form
\begin{equation}\label{e-1.1}
	h(z)=z+\sum_{n=2}^{\infty}a_nz^n\quad\mbox{and}\quad g(z)=\sum_{n=1}^{\infty}b_nz^n.
\end{equation}\vspace{1.5mm}

\par A function $ f\in\mathcal{H} $ is said to be convex, if $ f(\mathbb{D}) $ is convex. A domain $ \Omega $ is said to be close-to-convex if the complement of $ \Omega $ can be written as union of non-intersecting half lines. A function $ f\in\mathcal{H} $ is called close-to-convex in $ \mathbb{D} $, if $ f(\mathbb{D}) $ is close-to-convex. Mocanu \cite{Mocanu-1981} proved that if $ f $ is a harmonic mapping in a convex domain $ \Omega $ such that $ \Re(f_{z}(z))>|f_{\bar{z}}(z)| $ for $ z\in\Omega $, then the function $ f $ is not only univalent but also sense-preserving in $ \Omega $ while Ponnusamy \emph{et al.} \cite{Ponnusamy-CVEE-2013} showed that the functions in the class $ \mathcal{R}^{0}_{\mathcal{H}}=\{f=h+\bar{g}\in\mathcal{H} : \Re h^{\prime}(z)> |g^{\prime}(z)|\;\; \text{for all}\; z\in\mathbb{D}\} $ are close-to-convex in $ \mathbb{D} $.\vspace{1.2mm}

In the study of finding improved Bohr inequalities, the quantity $ S_r $ has some pivotal role. But in this paper, our prime concern is to study improved Bohr inequalities for certain classes of harmonic mappings using formula for $ S_r $. For analytic functions $ f $ defined on unit disk $ \mathbb{D} $, $S_r:=S_r(f)$ denotes the planar integral

$$S_r:=\int_{\mathbb D_r} |f'(z)|^2 d A(z),\;\;\mbox{where}\; 0<r<1.$$

If  $f(z)=\sum_{n=0}^\infty a_nz^n$, than the quantity $ S_r $ has the following series representation $ S_r:=\pi \sum_{n=1}^\infty n|a_n|^2 r^{2n}. $  In fact, if $f$ is a univalent function, then $S_r$ is the area of  $f(\mathbb D_r)$. In case of the multivalent function, $ S_r $ is greater than the area of the image of the subdisk $  \mathbb{D}_r $. This fact could be shown by noting that 
\begin{align*}
	S_r=\int_{\mathbb{D}_r}|f^{\prime}(w)|^2dA(w)=\int_{f(\mathbb{D}_r)}\nu_f(w)dA(w)\geq \int_{f(\mathbb{D}_r)}dA(w)=\mbox{Area}(f(\mathbb{D}_r)),
\end{align*}
where $ \nu_f(w)=\sum_{f(z)=w} 1$ denotes the counting function of $ f $.\vspace{1.2mm}

The quantity $ S_r $ plays a significant role in the study of improved versions of the classical Bohr inequality and with the help of this quantity, in the recent years, many  Bohr-type inequalities are obtained (see e.g. \cite{Ahamed-AASFM-2022,Ismagilov-2020-JMAA,Liu-Ponnusamy-BSM-2021,Huang-Liu-Ponnu-CVEE-2021}). For example, in $ 2020 $, Kayumov and Ponnusamy \cite{Kayumov-CRACAD-2020} established the following improved version of Bohr inequality for analytic functions.
\begin{thm}\cite{Kayumov-CRACAD-2020}\label{th-1.9}
	Let $ f(z)=\sum_{n=0}^{\infty}a_nz^n $ be analytic in $ \mathbb{D} $, $ |f(z)|\leq 1 $ and $ S_r $ denotes the image of the subdisk $ |z|<r $ under the mapping $ f $. Then 	$ B_1(r):=\sum_{n=0}^{\infty}|a_n|r^n+\frac{16}{9}\left(\frac{S_r}{\pi}\right)\leq 1\quad\mbox{for}\quad r\leq\frac{1}{3}, $
and the number $ 1/3 $ and $ 16/9 $ cannot be improved. Moreover, $ B_2(r):=|a_0|^2+\sum_{n=1}^{\infty}|a_n|r^n+\frac{8}{9}\left(\frac{S_r}{\pi}\right)\leq 1\quad\mbox{for}\quad r\leq\frac{1}{2}, $ 
and the number $ 1/2 $ and $ 8/9 $ cannot be improved.
\end{thm}
In $ 2020 $, Ismagilov \emph{et al.} \cite{Ismagilov-2020-JMAA} investigated further on Theorem \ref{th-1.9} and obtained the following result.
\begin{thm}\cite{Ismagilov-2020-JMAA}\label{th-1.10}
	Suppose that $ f(z)=\sum_{n=0}^{\infty}a_nz^n $ is analytic in $ \mathbb{D} $ and $ |f(z)|<1 $ in $ \mathbb{D} $. Then 
	\begin{equation*}
		M(r):=\sum_{n=0}^{\infty}|a_n||z|^n+\frac{16}{9}\left(\frac{S_r}{\pi}\right)+\lambda\left(\frac{S_r}{\pi}\right)^2\leq 1\;\; \mbox{for}\;\; r\leq\frac{1}{3}
	\end{equation*}
	where
	\begin{equation*}
		\lambda=\frac{4(486-261a-324a^2+2a^3+30a^4+3a^5)}{81(1+a)^3(3-5a)}=18.6095...
	\end{equation*}
	and $ a\approx 0.567284 $, is the unique positive root of the equation $ \psi(r)=0 $ in the interval $ (0,1) $, where $ \phi(t)=-405+473r+402r^2+38r^3+3r^4+r^5. $ The equality is achieved for the function $ f(z)=(a-z)/(1-az) $.
\end{thm}
It is natural to investigate on the following problem.
\begin{ques}\label{q-1.1}
	Can we prove the harmonic analogue of Theorems \ref{th-1.9} and \ref{th-1.10}?
\end{ques}
In this paper, we establish Bohr-type inequalities for polynomial expression of the quantity $ S_r/\pi $ for certain class of harmonic mappings. Surprisingly, we observed that like analytic case, the sharp coefficients in $ S_r/\pi $ and $ (S_r/\pi)^2 $ in Theorems \ref{th-1.9} and \ref{th-1.10} can be replaced by any non-negative real numbers in case of harmonic mappings.
\section{Improved Bohr inequalities for the class $ \mathcal{W}^{0}_{\mathcal{H}}(\alpha) $}
In $ 1977 $, the class $ \mathcal{W}(\alpha) $ was first introduced by Chichra \cite{Chichra-PAMS-1977}, which consisting of normalized analytic functions $ h $, satisfying the condition $ {\rm Re} \left(h^{\prime}(z)+\alpha zh^{\prime\prime}(z)\right)>0 $ for $ z\in\mathbb{D} $ and $ \alpha\geq 0 $. Furthermore, in \cite{Chichra-PAMS-1977} it has shown that functions in the class $ \mathcal{W}(\alpha) $ constitute a subclass of close-to-convex functions in $ \mathbb{D} $. In $ 2014 $, Nagpal and Ravichandran \cite{Nagpal-Ravinchandran-2014-JKMS} obtained coefficient bounds for functions in the class $ \mathcal{W}^{0}_{\mathcal{H}} $ which is defined by
\begin{equation*}
	\mathcal{W}^{0}_{\mathcal{H}}=\{f=h+\bar{g}\in \mathcal{H} :  {\rm Re}\left(h^{\prime}(z)+zh^{\prime\prime}(z)\right) > |g^{\prime}(z)+zg^{\prime\prime}(z)|\;\; \mbox{for}\; z\in\mathbb{D}\}.
\end{equation*}
It is known that $ \mathcal{W}^{0}_{\mathcal{H}} $ is closed under convolution and convex combinations. In $ 2019 $, Ghosh and Allu \cite{Nirupam-MonatsMath-2019} studied the class $ \mathcal{W}^{0}_{\mathcal{H}}(\alpha) $, where $ \alpha\geq 0 $ and 
\begin{equation*}
	\mathcal{W}^{0}_{\mathcal{H}}(\alpha)=\{f=h+\bar{g}\in \mathcal{H} :  {\rm Re}\left(h^{\prime}(z)+\alpha zh^{\prime\prime}(z)\right) > |g^{\prime}(z)+\alpha zg^{\prime\prime}(z)|\;\; \mbox{for}\; z\in\mathbb{D}\}
\end{equation*}
with $ g^{\prime}(0)=0 $ and proved that a harmonic mapping $ f=h+\bar{g} $ belongs to $ \mathcal{W}^{0}_{\mathcal{H}}(\alpha) $ if, and only if, the analytic function $ F=h+\epsilon g $ belongs to $ \mathcal{W}(\alpha) $ for each $ |\epsilon|=1. $ Clearly, $ \mathcal{W}^{0}_{\mathcal{H}}(\alpha) $ is a subclass of the close-to-convex harmonic  mappings. In particular, for $ \alpha\geq 1 $, the functions in the class $ \mathcal{W}^{0}_{\mathcal{H}} $ are fully starlike. In fact, the class  $ \mathcal{W}^{0}_{\mathcal{H}} $ is closed under convex combinations, and closed under convolutions. \vspace{1.2mm} 

\par Our purpose of this paper is to find sharp Bohr-type inequalities  for the class $ \mathcal{W}^{0}_{\mathcal{H}}(\alpha) $. Henceforth, we recall here the sharp coefficient bounds and the sharp growth estimates for functions in the class $ \mathcal{W}^{0}_{\mathcal{H}}(\alpha) $, which will play a key role to study the Bohr radius in terms of Euclidean distance.
\begin{lem}\cite{Nirupam-MonatsMath-2019}\label{lem-1.3}
	Let $ f\in \mathcal{W}^{0}_{\mathcal{H}}(\alpha) $ for $ \alpha\geq 0 $ and be of the form \eqref{e-1.1}. Then for any $ n\geq 2 $,
	\begin{enumerate}
		\item[(i)] $ |a_n|+|b_n|\leq \displaystyle\frac{2}{\alpha n^2+(1-\alpha)n} $;\vspace{1.5mm}
		\item[(ii)] $  ||a_n|-|b_n||\leq \displaystyle\frac{2}{\alpha n^2+(1-\alpha)n} $; \vspace{1.5mm}
		\item[(iii)] $ |a_n|\leq \displaystyle\frac{2}{\alpha n^2+(1-\alpha)n} $.
	\end{enumerate}
	All these inequalities are sharp for the function $ f=f_{\alpha} $ given by 
	\begin{equation}\label{e-1.4}
		f_{\alpha}(z)=z+\sum_{n=2}^{\infty}\frac{2z^n}{\alpha n^2+(1-\alpha)n}.
	\end{equation}
\end{lem}
\begin{lem}\cite{Nirupam-MonatsMath-2019}\label{lem-1.5}
	Let $ f\in \mathcal{W}^{0}_{\mathcal{H}}(\alpha) $ and be of the form \eqref{e-1.1} with $ 0<\alpha\leq 1 $. Then
	\begin{equation}\label{e-1.6}
		|z|+\sum_{n=2}^{\infty}\frac{2(-1)^{n-1}|z|^n}{\alpha n^2+(1-\alpha)n}\leq |f(z)|\leq |z|+\sum_{n=2}^{\infty}\frac{2|z|^n}{\alpha n^2+(1-\alpha)n}.
	\end{equation}
	Both the inequalities are sharp for the function  $ f=f_{\alpha} $ given by \eqref{e-1.4}.
\end{lem}
In view of Lemmas \ref{lem-1.3} and \ref{lem-1.5}, in $ 2020 $, Allu and Halder \cite{Himadri-Vasu-P1} studied the sharp Bohr radius for the class $ \mathcal{W}^{0}_{\mathcal{H}}(\alpha) $ and obtained the result.
\begin{thm}\cite{Himadri-Vasu-P1}\label{th-1.8}
	Let $ f\in \mathcal{W}^{0}_{\mathcal{H}}(\alpha) $ for $ \alpha\geq 0 $ be of the form \eqref{e-1.1}. Then 
	\begin{equation*}
		|z|+\sum_{n=2}^{\infty}\left(|a_n|+|b_n|\right)|z|^n\leq d(f(0),\partial f(\mathbb{D}))
	\end{equation*} 
	for $ |z|=r\leq r_f $, where $ r_f $ is the unique root of 
	\begin{equation*}
		r+\sum_{n=2}^{\infty}\frac{2r^n}{\alpha n^2+(1-\alpha)n}=1+\sum_{n=2}^{\infty}\frac{2(-1)^{n-1}}{\alpha n^2+(1-\alpha)n}
	\end{equation*}
	in $ (0,1) $. The radius $ r_f $ is the best possible.
\end{thm}
Before stating the main result of this paper, we first recall the definition of dilogarithm $ {\rm Li}_2(z) $ which is defined by the power series
\begin{equation*}
	{\rm Li}_2(z)=\sum_{n=1}^{\infty}\frac{z^n}{n^2}\;\; \mbox{for}\;\; |z|<1.
\end{equation*}
The definition and the name, of course, come from the analogy with the Taylor series of the ordinary logarithm around $ 1 $, 
\begin{equation*}
	-\log(1-z)=\sum_{n=1}^{\infty}\frac{z^n}{n}\;\;\mbox{for}\;\; |z|<1,
\end{equation*}
which leads similarly to the definition of the polylogarithm
\begin{equation*}
	{\rm Li}_m(z)=\sum_{n=1}^{\infty}\frac{z^n}{n^m}\;\; \mbox{for}\;\; |z|<1,\;\; m=1, 2, 3, \ldots.
\end{equation*}
In particular, the analytic continuation of the dilogarithm is given by
\begin{equation*}
	{\rm Li}_2(z)=-\int_{0}^{z}\log (1-u)\frac{du}{u}\;\; \mbox{for}\;\; z\in\mathbb{C}\setminus [1,\infty).
\end{equation*}
The quantity $ S_r $ has been defined for certain class of harmonic mappings (see \cite{Aha-CMFT-2022,Ahamed-AMP-2021,Ahamed-CVEE-2021,Allu-Hal-IM-2021,Wang-AMP-2022}) and improved Bohr radii are studied for such classes. In this paper, we aim to study improved Bohr radius for the class $ \mathcal{W}^{0}_{\mathcal{H}}(\alpha) $ and to show that more general improved versions of the Bohr inequality can be established. Before we state the main result, we define here the polynomial $ P $ of degree $ k $ as the following
\begin{align*}
	P(w)=\lambda_kw^k+\cdots+\lambda_1w,\; \mbox{where}\; \lambda_k(\neq 0),\;\lambda_j\in\{x\in\mathbb{R}: x\geq 0\}.
\end{align*}
 Using the Lemmas \ref{lem-1.3} and \ref{lem-1.5}, we obtain the following sharp Bohr-type inequalities for the class $\mathcal{W}^{0}_{\mathcal{H}}(\alpha) $.
 \begin{thm}\label{th-2.13}
 	Let $ f\in \mathcal{W}^{0}_{\mathcal{H}}(\alpha)  $, where $ 0<\alpha\leq 1 $, be given by \eqref{e-1.1}. Then 
 	\begin{enumerate}
 		\item[(i)]  the inequality $ 	|z|+\sum_{n=2}^{\infty}\left(|a_n|+|b_n|\right)|z|^n+P\left(\frac{S_r}{\pi}\right)\leq d(f(0),\partial f(\mathbb{D})) $
 		holds for $ |z|=r\leq r_k(\alpha) $, where $ r_k(\alpha) $ is the unique root of the equation $ J_1(r)=0 $ and
 		\begin{align*}
 			J_1(r):&=r+\sum_{n=2}^{\infty}\frac{2r^n}{\alpha n^2+(1-\alpha)n}+P\left(r^2+\sum_{n=2}^{\infty}\frac{4nr^{2n}}{(\alpha n^2+(1-\alpha)n)^2}\right)-1\\&\quad\quad-\sum_{n=2}^{\infty}\frac{2(-1)^{n-1}}{\alpha n^2+(1-\alpha)n}=0
 		\end{align*}
 		in $ (0,1) $. Here $ r_k(\alpha) $ is the best possible.\vspace{1.5mm}
 		\item[(ii)] for integer $ m\geq 0 $, the inequality $ 	|f(z)|^m+\sum_{n=2}^{\infty}\left(|a_n|+|b_n|\right)|z|^n+P\left(\frac{S_r}{\pi}\right)\leq d(f(0),\partial f(\mathbb{D})) $
 		holds for $ |z|=r\leq r^{*}_k(m, \alpha) $, where $ r^{*}_k(m,\alpha) $ is the unique root of the equation $ J_2(r)=0 $ and
 		\begin{align*}
 			J_2(r):=&\left(r+\sum_{n=2}^{\infty}\frac{2r^n}{\alpha n^2+(1-\alpha)n}\right)^m+P\left(r^2+\sum_{n=2}^{\infty}\frac{4nr^{2n}}{(\alpha n^2+(1-\alpha)n)^2}\right)\\&\quad\quad+\sum_{n=2}^{\infty}\frac{2r^n}{\alpha n^2+(1-\alpha)n}-1-\sum_{n=2}^{\infty}\frac{2(-1)^{n-1}}{\alpha n^2+(1-\alpha)n}=0
 		\end{align*}
 		in $ (0,1) $.  Here $ r^{*}_k(m, \alpha) $ is the best possible.
 	\end{enumerate}
 \end{thm} 
\begin{cor}\label{cor-2.7}
	Let $ f\in \mathcal{W}^{0}_{\mathcal{H}}(\alpha)  $, where $ 0<\alpha\leq 1 $, be given by \eqref{e-1.1}. Then   the inequality $ 	|z|+\sum_{n=2}^{\infty}\left(|a_n|+|b_n|\right)|z|^n+\lambda_1\left(\frac{S_r}{\pi}\right)+\lambda_2\left(\frac{S_r}{\pi}\right)^2\leq d(f(0),\partial f(\mathbb{D})) $ holds for $ |z|=r\leq r_2(\alpha) $, where $ r_2(\alpha) $ is the unique root of the equation $ J_1(r)=0 $ and
		\begin{align*}
			&r+\sum_{n=2}^{\infty}\frac{2r^n}{\alpha n^2+(1-\alpha)n}+\lambda_1\left(r^2+\sum_{n=2}^{\infty}\frac{4nr^{2n}}{(\alpha n^2+(1-\alpha)n)^2}\right)\\&\quad+\lambda_2\left(r^2+\sum_{n=2}^{\infty}\frac{4nr^{2n}}{(\alpha n^2+(1-\alpha)n)^2}\right)^2=1+\sum_{n=2}^{\infty}\frac{2(-1)^{n-1}}{\alpha n^2+(1-\alpha)n}=0
		\end{align*}
		in $ (0,1) $. Here $ r_2(\alpha) $ is the best possible.
\end{cor}
\begin{rem}
	For the class $ \mathcal{W}^{0}_{\mathcal{H}}(\alpha) $, the Corollary \ref{cor-2.7} is the harmonic analogue of Theorem \ref{th-1.10} for $ \lambda_1=16/9 $ and $ \lambda_2=\lambda $ as in Theorem \ref{th-1.10}. Thus Question \ref{q-1.1} is answered completely for the class $ \mathcal{W}^{0}_{\mathcal{H}}(\alpha) $.
\end{rem}
As a corollary of Theorem \ref{th-2.13}, we have the following result.
\begin{cor}
	Let $ f\in \mathcal{W}^{0}_{\mathcal{H}}(\alpha)  $, where $ 0<\alpha\leq 1 $, be given by \eqref{e-1.1}. Then 
	\begin{enumerate}
		\item[(i)]  the inequality $ 	|z|+\sum_{n=2}^{\infty}\left(|a_n|+|b_n|\right)|z|^n+\frac{S_r}{\pi}\leq d(f(0),\partial f(\mathbb{D})) $
	holds for $ |z|=r\leq r_f(\alpha) $, where $ r_f(\alpha) $ is the unique root of the equation 
\begin{align*}
	&r^2+r+\sum_{n=2}^{\infty}\frac{2r^n}{\alpha n^2+(1-\alpha)n}+\sum_{n=2}^{\infty}\frac{4nr^{2n}}{(\alpha n^2+(1-\alpha)n)^2}-1\\&\quad\quad-\sum_{n=2}^{\infty}\frac{2(-1)^{n-1}}{\alpha n^2+(1-\alpha)n}=0
\end{align*}
in $ (0,1) $. The constant $ r_f(\alpha) $ is the best possible.\vspace{1.5mm}
\item[(ii)] for integer $ m\geq 0 $, the inequality $ |f(z)|^2+\sum_{n=2}^{\infty}\left(|a_n|+|b_n|\right)|z|^n+\left(\frac{S_r}{\pi}\right)^k\leq d(f(0),\partial f(\mathbb{D})) $ holds for $ |z|=r\leq r^{*}_f(k, \alpha) $, where $ r^{*}_f(k,\alpha) $ is the unique root of the equation 
\begin{align*}
	&\left(r+\sum_{n=2}^{\infty}\frac{2r^n}{\alpha n^2+(1-\alpha)n}\right)^2+\left(r^2+\sum_{n=2}^{\infty}\frac{4nr^{2n}}{(\alpha n^2+(1-\alpha)n)^2}\right)^m\\&\quad\quad+\sum_{n=2}^{\infty}\frac{2r^n}{\alpha n^2+(1-\alpha)n}-1-\sum_{n=2}^{\infty}\frac{2(-1)^{n-1}}{\alpha n^2+(1-\alpha)n}=0
\end{align*}
in $ (0,1) $.  The constant $ r^{*}_f(k, \alpha) $ is best possible.
	\end{enumerate}
\end{cor} 
As a corollary of Theorem \ref{th-2.13}, we obtain the following result finding sharp Bohr radius for the $ \mathcal{W}^{0}_{\mathcal{H}}(\alpha) $.
\begin{cor}\label{cor-2.14}
		Let $ f\in \mathcal{W}^{0}_{\mathcal{H}}(\alpha)  $ be given by \eqref{e-1.1}. Then
	\begin{enumerate}
		\item[(i)]   for $ \alpha=1/2 $, the inequality $ |z|+\sum_{n=2}^{\infty}\left(|a_n|+|b_n|\right)|z|^n+\frac{S_r}{\pi}\leq d(f(0),\partial f(\mathbb{D})) $
		holds for $ |z|=r\leq r_f(\alpha)\approx 0.600881 $, where $ r_f(\alpha) $ is the unique root of the equation $ F(r)=0 $ and 
		\begin{align*}
		F(r):=&\frac{4}{r}(1-r)\log (1-r)+\frac{16}{r^2}(1-r^2)\log(1-r^2)-\frac{16}{r^2}{\rm Li}_2(r^2)\\&\quad\quad-3r^2-r+29+8\log 2=0
		\end{align*}
		in $ (0,1) $. Here $ r_f(\alpha)\approx 0.600881 $ is the best possible.\vspace{1.5mm}
		\item[(ii)]  for $ \alpha=1/2 $, the inequality $ |f(z)|+\sum_{n=2}^{\infty}\left(|a_n|+|b_n|\right)|z|^n+\frac{S_r}{\pi}\leq d(f(0),\partial f(\mathbb{D})) $
		holds for $ |z|=r\leq r^{*}_f(\alpha)\approx 0.302059 $, where $ r^{*}_f(\alpha) $ is the unique root of the equation $ T(r)=0 $ and 
		\begin{align*}
			T(r):=&\frac{8}{r}(1-r)\log (1-r)+\frac{16}{r^2}(1-r^2)\log (1-r^2)-\frac{16}{r^2}{\rm Li}_2(r^2)\\&\quad\quad -3r^2-3r+45-8\log 2=0
		\end{align*}
		in $ (0,1) $. Here $ r^{*}_f(\alpha)\approx 0.302059 $ is the best possible.
	\end{enumerate}
\end{cor}
\section{Improved Bohr inequalities for stable harmonic mappings}
A log-harmonic mapping defined on the unit disk $ \mathbb{D} $ is a solution of the non-linear elliptic partial differential equation
\begin{align}
	\frac{\overline{f_{\bar{z}}(z)}}{\overline{f(z)}}=\omega(z)\frac{f_{z}(z)}{f(z)},
\end{align}
where $ \omega $ is an analytic function defined on $ \mathbb{D} $. It is well-known that if $ f $ is non-constant log-harmonic mapping of $ \mathbb{D} $ and vanishes only at $ z=0 $, then $ f $ admits the representation 
\begin{align*}
	f(z)=z^m|z|^{2\beta m}h(z)\overline{g(z)},
\end{align*}
where $ m $ is a non-negative integer, $ Re(\beta)>-\frac{1}{2} $, $ h $ and $ g $ are analytic functions in $ \mathbb{D} $ satisfying $ h(0)\neq 0 $ and $ g(0)=1 $. For more detailed, we refer to the article \cite{Abdulhadi-TAMS-1988}. It is worth pointing out that $ f(0)\neq 0 $ if, and only if, $ m=1 $, that is, $ f $ has the following representation 
\begin{align*}
	f(z)=z|z|^{2\beta}h(z)\overline{g(z)},
\end{align*}
where $ Re(\beta)>-\frac{1}{2} $ and $ 0\not\in (hg)(\mathbb{D}) $.\vspace{1.2mm}

We now recall here the definition of stable harmonic univalent and stable convex harmonic mappings.
\begin{defn}\cite{Hernandez-PCPS-2013}
	A (sense preserving) harmonic mapping $ f=h+\bar{g} $ is said to be stable harmonic univalent or $ SHU $ in the unit disk $ \mathbb{D} $ (resp. stable convex harmonic or $ SHC $) if all the mappings $ f_{\lambda}=h+\lambda g $ with $ |\lambda|=1 $ are univalent (resp. convex) in $ \mathbb{D} $.
\end{defn}
The class $ \mathcal{S}^{0}_{\mathcal{H}} $ is defined as the family of sense preserving univalent harmonic mappings $ f=h+\bar{g} $ in the unit disk $ \mathbb{D} $ with the normalization $ h(0)=g(0)=1-h^{\prime}(0)=g^{\prime}(0)=0 $. Hence, we suppose that $ 	h(z)=z+\sum_{n=2}^{\infty}a_nz^n\;\mbox{and}\; g(z)=\sum_{n=2}^{\infty}b_nz^n. $ The majorant series corresponding to $ f=h+\bar{g} $ is defined as 
\begin{align*}
	M_f(r)=r+\sum_{n=2}^{\infty}\left(|a_n|+|b_n|\right)r^n,\; \mbox{where}\; |z|=r .
\end{align*}
The key ingredient in our study is the coefficient bounds and distortion theorems which were proved by in \cite{Hernandez-PCPS-2013}.
\begin{lem}\cite{Hernandez-PCPS-2013}\label{lem-3.2}
	\begin{enumerate}
		\item[(i)] Assume that $ f=h+\bar{g} $ in $ \mathcal{S}^{0}_{\mathcal{H}} $
		is stable convex harmonic mapping. Then for all non-negative integers $ n $, we have
		\begin{align*}
			||a_n|-|b_n||\leq \max\{|a_n|, |b_n|\}\leq |a_n|+|b_n|\leq 1.
		\end{align*}
		\item[(ii)] Assume that $ f=h+\bar{g} $ in $ \mathcal{S}^{0}_{\mathcal{H}} $
		is stable univalent harmonic mapping. Then for all non-negative integers $ n $, we have
		\begin{align*}
			||a_n|-|b_n||\leq \max\{|a_n|, |b_n|\}\leq |a_n|+|b_n|\leq n.
		\end{align*} 
		\item[(iii)] If $ f=h+\bar{g} $ in $ \mathcal{S}^{0}_{\mathcal{H}} $
		is a stable convex harmonic mapping in the unit disk $ \mathbb{D} $, then we have
		\begin{align}\label{e-3.4}
			\frac{|z|}{1+|z|}\leq |f(z)|\leq\frac{|z|}{1-|z|}.
		\end{align}
	\item[(iv)] If $ f=h+\bar{g} $ in $ \mathcal{S}^{0}_{\mathcal{H}} $
	is a stable univalent harmonic mapping in the unit disk $ \mathbb{D} $, then we have
	\begin{align}\label{e-3.3}
		\frac{|z|}{(1+|z|)^2}\leq |f(z)|\leq\frac{|z|}{(1-|z|)^2}.
	\end{align}
	\end{enumerate}

\end{lem}
In view of Lemma \ref{lem-3.2}, a simple computation shows that when $ f $ is a stable univalent mapping, then $ S_r/\pi\leq \sum_{n=1}^{\infty}nr^{2n}=r^2/(1-r^2)^2 $ and if $ f $ is stable convex harmonic mapping, then $ S_r/\pi\leq \sum_{n=1}^{\infty}n^3r^{2n}=(r^6+4r^4+r^2)/(r^2-1)^2 $. In view of the bound of the quantity $ S_r/\pi $, 
Abdulhadi and Hajj \cite{Abdulhadi-Hajj-2022}, recently, proved a version of the Improved Bohr’s inequality under the stability condition for the harmonic mappings by adding a suitable non-negative term at the left hand side of the inequality.
\begin{thm}\cite{Abdulhadi-Hajj-2022}\label{th-3.7}
	\begin{enumerate}
		\item[(i)] Let $ f=h+\bar{g}\in\mathcal{S}^{0}_{\mathcal{H}} $ be a stable convex harmonic mapping on the unit disk $ \mathbb{D} $, and let $ S_r $ be the area of the image $ f(\mathbb{D}_r) $, with $ \mathbb{D}_r $ is a sub-disk of $ \mathbb{D} $. Then
		\begin{align*}
			M_f(r)+\left(\frac{S_r}{\pi}\right)^k\leq d(f(0), \partial f(\mathbb{D}))
		\end{align*}
		for $ |z|\leq r_0 $, where $ r_0 $ is the unique root in $ (0, 1) $ of 
		\begin{align*}
			\frac{r}{1-r}+\frac{r^{2k}}{(1-r^2)^{2k}}=\frac{1}{2}.
		\end{align*}
		Here $ r_0 $ is best possible.
		\item[(ii)] Let $ f=h+\bar{g}\in\mathcal{S}^{0}_{\mathcal{H}} $ be a stable univalent harmonic mapping on the unit disk $ \mathbb{D} $. Then
		\begin{align*}
			M_f(r)+\left(\frac{S_r}{\pi}\right)^k\leq d(f(0), \partial f(\mathbb{D}))
		\end{align*}
		for $ |z|\leq r_0 $, where $ r_0\approx 0.382 $ is the unique root in $ (0, 1) $ of 
		\begin{align*}
			\frac{r}{(1-r)^2}+\frac{\left(r^6+4r^4+r^2\right)^k}{(r^2-1)^{4k}}=\frac{1}{4}.
		\end{align*}
		Here $ r_0 $ is best possible.
	\end{enumerate}
\end{thm}
In $ 2021 $, strengthening Theorem \ref{th-1.9}, Ismagilov \emph{et al.} \cite{Ismagilov-JMS-2021} proved the following result.
\begin{thm}\cite{Ismagilov-JMS-2021}\label{th-3.6}
	Let $ f(z)=\sum_{n=0}^{\infty}a_nz^n $ be analytic in $ \mathbb{D} $, $ |f(z)|\leq 1 $ and $ S_r $ denotes the image of the sub-disk $ |z|<r $ under the mapping $ f $. Then 
	\begin{align*}
		\sum_{n=0}^{\infty}|a_n|r^n+\frac{16}{9}\left(\frac{S_r}{\pi-S_r}\right)\leq 1\; \mbox{for}\; r\leq \frac{1}{3},
	\end{align*}
	and the number $ 16/9 $ cannot be improved. 
\end{thm}
In this section, we have two objectives: firstly, we generalize Theorem \ref{th-3.7} as an improved version of Bohr inequalities in terms adding a polynomial expression $ P(S_r/\pi) $ with the majorant series $ M_f(r) $, and secondly we give a version of Theorem \ref{th-3.7} terms adding a polynomial expression $ P(S_r/(\pi-S_r)) $ with the majorant series $ M_f(r) $.\vspace{1.2mm}

\noindent We obtain the following result generalizing  Theorem \ref{th-3.7}.
\begin{thm}\label{th-3.8}
	\begin{enumerate}
		\item[(i)] Let $ f=h+\bar{g}\in\mathcal{S}^{0}_{\mathcal{H}} $ be a stable convex harmonic mapping on the unit disk $ \mathbb{D} $, and let $ S_r $ be the area of the image $ f(\mathbb{D}_r) $, with $ \mathbb{D}_r $ is a sub-disk of $ \mathbb{D} $. Then
		\begin{align*}
			M_f(r)+P\left(\frac{S_r}{\pi}\right)\leq d(f(0), \partial f(\mathbb{D}))
		\end{align*}
		for $ |z|\leq r_k $, where $ r_k $ is the unique root in $ (0, 1) $ of 
		\begin{align*}
			\frac{r}{1-r}+\frac{\lambda_k r^{2k}}{(1-r^2)^{2k}}+\cdots+\frac{\lambda_1 r^{2}}{(1-r^2)^{2}}=\frac{1}{2}.
		\end{align*}
		Here $ r_k $ is best possible.
		\item[(ii)] Let $ f=h+\bar{g}\in\mathcal{S}^{0}_{\mathcal{H}} $ be a stable univalent harmonic mapping on the unit disk $ \mathbb{D} $. Then
		\begin{align*}
			M_f(r)+P\left(\frac{S_r}{\pi}\right)\leq d(f(0), \partial f(\mathbb{D}))
		\end{align*}
		for $ |z|\leq r^{*}_k $, where $ r^{*}_k $ is the unique root in $ (0, 1) $ of 
		\begin{align*}
			\frac{r}{(1-r)^2}+\frac{\lambda_k\left(r^6+4r^4+r^2\right)^k}{(r^2-1)^{4k}}+\cdots+\frac{\lambda_1\left(r^6+4r^4+r^2\right)}{(r^2-1)^{4}}=\frac{1}{4}.
		\end{align*}
		Here $ r^{*}_k $ is best possible.
	\end{enumerate}
\end{thm}
\begin{rem}
	It is easy to see that, in particular when $ P(z)=z^k $, then Theorem \ref{th-3.7} is a special case of Theorem \ref{th-3.8}. Thus, Theorem \ref{th-3.8} generalizes Theorem \ref{th-3.7}.
\end{rem}
\begin{cor}\label{cor-3.8}
	\begin{enumerate}
		\item[(i)] Let $ f=h+\bar{g}\in\mathcal{S}^{0}_{\mathcal{H}} $ be a stable convex harmonic mapping on the unit disk $ \mathbb{D} $, and let $ S_r $ be the area of the image $ f(\mathbb{D}_r) $, with $ \mathbb{D}_r $ is a sub-disk of $ \mathbb{D} $. Then
		\begin{align*}
			M_f(r)+\lambda_1\left(\frac{S_r}{\pi}\right)+\lambda_2\left(\frac{S_r}{\pi}\right)^2\leq d(f(0), \partial f(\mathbb{D}))
		\end{align*}
		for $ |z|\leq r_k $, where $ r_k $ is the unique root in $ (0, 1) $ of 
		\begin{align*}
			\frac{r}{1-r}+\frac{\lambda_2 r^{4}}{(1-r^2)^{4}}+\frac{\lambda_1 r^{2}}{(1-r^2)^{2}}=\frac{1}{2}.
		\end{align*}
		Here $ r_k $ is best possible.
		\item[(ii)] Let $ f=h+\bar{g}\in\mathcal{S}^{0}_{\mathcal{H}} $ be a stable univalent harmonic mapping on the unit disk $ \mathbb{D} $. Then
		\begin{align*}
			M_f(r)+\lambda_1\left(\frac{S_r}{\pi}\right)+\lambda_2\left(\frac{S_r}{\pi}\right)^2\leq d(f(0), \partial f(\mathbb{D}))
		\end{align*}
		for $ |z|\leq r^{*}_k $, where $ r^{*}_k $ is the unique root in $ (0, 1) $ of 
		\begin{align*}
			\frac{r}{(1-r)^2}+\frac{\lambda_2\left(r^6+4r^4+r^2\right)^2}{(r^2-1)^{8}}+\frac{\lambda_1\left(r^6+4r^4+r^2\right)}{(r^2-1)^{4}}=\frac{1}{4}.
		\end{align*}
		Here $ r^{*}_k $ is best possible.
	\end{enumerate}
\end{cor}
\begin{rem}
	We see that Corollary \ref{cor-3.8} is the harmonic analogue of Theorem \ref{th-1.10} for the class $ \mathcal{S}^{0}_{\mathcal{H}} $ for $ \lambda_1=16/9 $ and $ \lambda_2=\lambda $ as in Theorem \ref{th-1.10}. Thus the Question \ref{q-1.1} is answered for the class $ \mathcal{S}^{0}_{\mathcal{H}} $.
\end{rem}
To the best of our knowledge, no study on Bohr inequalities are done yet for the class of harmonic mappings. In this section, our purpose is to study improved Bohr inequalities for the class of stable harmonic mappings. Hence, motivated from the result of Ismagilov \emph{et al.} \cite{Ismagilov-JMS-2021}, we see that 
 $ \frac{S_r}{\pi}\leq \frac{r^2}{(1-r^2)^2} $ for $ r\in (0, 1) $ for the class $ \mathcal{S}^{0}_{\mathcal{H}} $ when $ f $ is stable convex harmonic, we see that $ \frac{\pi}{\pi-S_r}\leq \frac{(1-r^2)^2}{(1-r^2)^2-r^2} $. Therefore, it follows that 
 \begin{align}\label{e-33.8}
 	\dfrac{S_r}{\pi-S_r}\leq\dfrac{r^2}{(1-r^2)^2-r^2}\; \mbox{for} r\in (0, 1).
 \end{align} Similarly, when $ f\in \mathcal{S}^{0}_{\mathcal{H}} $ is stable univalent harmonic, it can be shown that 
\begin{align}\label{e-33.9}
	\dfrac{S_r}{\pi-S_r}\leq \dfrac{r^6+4r^4+r^2}{(r^2-1)^4-(r^6+4r^4+r^2)} \; \mbox{for} r\in (0, 1).
\end{align}\vspace{1.2mm} 
 
 We obtain the next result using the upper bounds of the quantity $ \dfrac{S_r}{\pi-S_r} $ for functions in the class $ \mathcal{S}^{0}_{\mathcal{H}} $.
\begin{thm}\label{th-3.9}
	\begin{enumerate}
		\item[(i)] Let $ f=h+\bar{g}\in\mathcal{S}^{0}_{\mathcal{H}} $ be a stable convex harmonic mapping on the unit disk $ \mathbb{D} $, and let $ S_r $ be the area of the image $ f(\mathbb{D}_r) $, with $ \mathbb{D}_r $ is a sub-disk of $ \mathbb{D} $. Then
		\begin{align*}
			M_f(r)+P\left(\frac{S_r}{\pi-S_r}\right)\leq d(f(0), \partial f(\mathbb{D}))
		\end{align*}
		for $ |z|\leq R_k $, where $ R_k $ is the unique root in $ (0, 1) $ of 
		\begin{align*}
			\frac{r}{1-r}+\dfrac{\lambda_k r^{2k}}{((1-r^2)^2-r^2)^k}+\cdots+\dfrac{\lambda_1r^2}{(1-r^2)^2-r^2}=\frac{1}{2}.
		\end{align*}
		Here $ R_k $ is best possible.
		\item[(ii)] Let $ f=h+\bar{g}\in\mathcal{S}^{0}_{\mathcal{H}} $ be a stable univalent harmonic mapping on the unit disk $ \mathbb{D} $. Then
		\begin{align*}
			M_f(r)+P\left(\frac{S_r}{\pi-S_r}\right)\leq d(f(0), \partial f(\mathbb{D}))
		\end{align*}
		for $ |z|\leq R^{*}_k $, where $ R_k $ is the unique root in $ (0, 1) $ of 
		\begin{align*}
			&\frac{r}{(1-r)^2}+\lambda_k\left(\dfrac{r^6+4r^4+r^2}{(r^2-1)^4-(r^6+4r^4+r^2)}\right)^k\\&+\cdots+\lambda_1\left(\dfrac{r^6+4r^4+r^2}{(r^2-1)^4-(r^6+4r^4+r^2)}\right)=\frac{1}{4}.
		\end{align*}
		Here $ R^{*}_k $ is best possible.
	\end{enumerate}
\end{thm}
\section{Proof of the main results}
\begin{proof}[\bf Proof of Theorem \ref{th-2.13}]
	For $ f\in\mathcal{W}^0_{\mathcal{H}}(\alpha) $, the Jacobian of $ f $ is denoted by $ J_f $ and is defined by
	\begin{equation*}
		J_f(z)=|f_{z}(z)|^2-|f_{\bar{z}}(z)|^2=|h^{\prime}(z)|^2-|g^{\prime}(z)|^2\;\; \mbox{for}\;\; z\in\mathbb{D}.
	\end{equation*}
It is well-known that (see \cite[p.113]{Duren-Harmonic-2004} and see also \cite{Evdoridis-IM-2018}) the area of the disk $ \mathbb{D}_r:=\{z\in\mathbb{C} : |z|<r\} $ under the harmonic map $ f=h+\bar{g} $ is defined by
\begin{align}\label{e-3.18}
	S_r=\iint\limits_{\mathbb{D}_r}J_f(z)dxdy=\iint\limits_{\mathbb{D}_r}\left(|{h^{\prime}(z)}|^2-|{g^{\prime}(z)}|^2\right)dxdy.
\end{align}
In addition, we have the following (see \cite[Proof of Theorem 3.2]{Ahamed-CVEE-2021} for detailed information)
\begin{align}\label{e-3.21}
	\frac{S_r}{\pi}=r^2+\sum_{n=2}^{\infty}\frac{4nr^{2n}}{\left(\alpha n^2+(1-\alpha)n\right)^2}.
\end{align}
Let $ f\in \mathcal{W}^{0}_{\mathcal{H}}(\alpha)  $ be given by \eqref{e-1.1}. Then, in view of Lemma \ref{lem-1.5} for $ |z|=r $, we have 
\begin{equation}\label{e-33.11}
	r+\sum_{n=2}^{\infty}\frac{2(-1)^{n-1}r^n}{\alpha n^2+(1-\alpha)n}\leq |f(z)|\leq r+\sum_{n=2}^{\infty}\frac{2r^n}{\alpha n^2+(1-\alpha)n}.
\end{equation}
It is easy to see that $ f(0)=0 $ and hence, $ |f(z)-f(0)|=|f(z)| $. Therefore, we have the following inequality
\begin{align}\label{ee-33.66}
	\liminf_{r\rightarrow 1^{-}}\left(r+\sum_{n=2}^{\infty}\frac{2(-1)^{n-1}r^n}{\alpha n^2+(1-\alpha)n}\right)	&\leq \liminf_{r\rightarrow 1^{-}}|f(z)-f(0)|\\&\nonumber\leq \liminf_{r\rightarrow 1^{-}}\left(r+\sum_{n=2}^{\infty}\frac{2r^n}{\alpha n^2+(1-\alpha)n}\right).
\end{align}
For $ 0<\alpha\leq 1 $ and $ n\geq 2 $, since $ \alpha n^2+(1-\alpha)n\geq \alpha n^2>0 $, hence, we obtain
\begin{equation*}
	\sum_{n=2}^{\infty}\frac{2r^n}{\alpha n^2+(1-\alpha)n}\leq \sum_{n=2}^{\infty}\frac{2r^n}{\alpha n^2}.
\end{equation*} 
We define the functions $ g_n $ and $ M_n $ by
\begin{equation*}
	g_n(r):=\frac{2r^n}{\alpha n^2+(1-\alpha)n}\;\; \mbox{and}\;\; M_n:=\frac{2r^n}{\alpha n^2}.
\end{equation*} 
Evidently, $ M_n>0 $ and $ |g_n(r)|\leq M_n $ for each $ n\geq 2 $. Since $ |z|=r<1 $, a simple computation shows that
\begin{equation*}
	\sum_{n=2}^{\infty}M_n<\frac{2}{\alpha }\sum_{n=2}^{\infty}\frac{1}{n^2}=\frac{2}{\alpha }\left(\frac{\pi^2}{6}-1\right).
\end{equation*}
By the comparison test, the series $ \sum_{n=2}^{\infty}M_n $ converges,  hence, by the Weierstrass $ M $-test for series of functions, the series $$ \sum_{n=2}^{\infty}g_n(r)=\sum_{n=2}^{\infty}\frac{2r^n}{\alpha n^2+(1-\alpha)n} $$ is absolutely and uniformly convergent in $ |z|=r<1 $ and hence, we interchange limit and the summation
\begin{equation*}
	\lim_{r\rightarrow 1^{-}}\sum_{n=2}^{\infty}g_n(r)=\sum_{n=2}^{\infty}\lim_{r\rightarrow 1^{-}}g_n(r).
\end{equation*}
Therefore, we must have
\begin{align*}
	\liminf_{r\rightarrow 1^{-}}\left(r+\sum_{n=2}^{\infty}\frac{2r^n}{\alpha n^2+(1-\alpha)n}\right)=1+\sum_{n=2}^{\infty}\frac{2}{\alpha n^2+(1-\alpha)n}.
\end{align*}

On the other hand, we see that
\begin{equation*}
	r+\sum_{n=2}^{\infty}\bigg|\frac{2(-1)^{n-1}r^n}{\alpha n^2+(1-\alpha)n}\bigg|\leq r+\sum_{n=2}^{\infty}\frac{2r^n}{\alpha n^2}=r+\sum_{n=2}^{\infty}M_n<1+\frac{2}{\alpha}\left(\frac{\pi^2}{6}-1\right),
\end{equation*}
and hence, by Weierstrass M-test, the series $ 	r+\sum_{n=2}^{\infty}\frac{2(-1)^{n-1}r^n}{\alpha n^2+(1-\alpha)n} $ is also absolutely and uniformly convergent in $ |z|=r<1 $. By using the same argument used in the above, we obtain
\begin{align}\label{eee-3.2}
	\liminf_{r\rightarrow 1^{-}}\left(r+\sum_{n=2}^{\infty}\frac{2(-1)^{n-1}r^n}{\alpha n^2+(1-\alpha)n}\right)=1+\sum_{n=2}^{\infty}\frac{2(-1)^{n-1}}{\alpha n^2+(1-\alpha)n}.
\end{align}
Thus, it follows from \eqref{ee-33.66} that
\begin{equation}\label{ee-33.33}
	1+\sum_{n=2}^{\infty}\frac{2(-1)^{n-1}}{\alpha n^2+(1-\alpha)n}\leq \liminf_{r\rightarrow 1^{-}}|f(z)-f(0)|\leq 1+\sum_{n=2}^{\infty}\frac{2}{\alpha n^2+(1-\alpha)n}.
\end{equation}
Therefore, the Euclidean distance $ d $ between $ f(0) $ and the boundary of $ f(\mathbb{D}) $ is 
\begin{align}
	\label{e-3.1}
	d(f(0),\partial f(\mathbb{D}))=\liminf_{|z|=r\rightarrow 1^{-}}|f(z)-f(0)|\geq 1+\sum_{n=2}^{\infty}\frac{2(-1)^{n-1}}{\alpha n^2+(1-\alpha)n}.
\end{align}
(i) In view of Lemmas \ref{lem-1.3} and \ref{lem-1.5}, and  \eqref{e-3.21} for $ |z|=r $, we obtain
\begin{align}\label{e-3.22}
	&|z|+\sum_{n=2}^{\infty}\left(|a_n|+|b_n|\right)|z|^n+P\left(\frac{S_r}{\pi}\right)\\&\nonumber\leq r+\sum_{n=2}^{\infty}\frac{2r^n}{\alpha n^2+(1-\alpha)n}+P\left(r^2+\sum_{n=2}^{\infty}\frac{4nr^{2n}}{\left(\alpha n^2+(1-\alpha)n\right)^2}\right).
\end{align}
It is easy to see that
\begin{align*}
	 r+\sum_{n=2}^{\infty}\frac{2r^n}{\alpha n^2+(1-\alpha)n}+P\left(r^2+\sum_{n=2}^{\infty}\frac{4nr^{2n}}{\left(\alpha n^2+(1-\alpha)n\right)^2}\right)\leq 1+\sum_{n=2}^{\infty}\frac{2(-1)^{n-1}}{\alpha n^2+(1-\alpha)n}
\end{align*}
for $ r\leq  r_k(\alpha)$, where $ r_k(\alpha) $ a the root of $ J_1(r)=0 $ in $ (0,1) $.\vspace{1.2mm}
\par It is not difficult to show that $ J_1(0)J_1(1)<0 $ and $ J^{\prime}_1(r)>0 $ for $ r\in (0,1) $. Then by the Intermediate value theorem, $ r_k(\alpha) $ is the unique roof of $ J_1(r)=0 $ in $ (0,1) $. Therefore, we have
\begin{align}\label{e-3.23}
	&r_f(\alpha)+\sum_{n=2}^{\infty}\frac{2r^n_f(\alpha)}{\alpha n^2+(1-\alpha)n}+P\left(r^2_f(\alpha)+\sum_{n=2}^{\infty}\frac{4nr^{2n}_f(\alpha)}{\left(\alpha n^2+(1-\alpha)n\right)^2}\right)\\&\nonumber= 1+\sum_{n=2}^{\infty}\frac{2(-1)^{n-1}}{\alpha n^2+(1-\alpha)n}.
\end{align}
\par In order to show that $ r_f(\alpha) $ is the best possible, we consider the function $ f=f_{\alpha} $ which is defined by \eqref{e-1.4}. Note that $ f_{\alpha}\in \mathcal{W}^{0}_{\mathcal{H}}(\alpha) $ and $ f_{\alpha}(0)=0 $. Now at $ z=-r $, we verify that
\begin{align}\label{ee-33.88}
	|f_{\alpha}(-r)-f_{\alpha}(0)|=\bigg|-r+\sum_{n=2}^{\infty}\frac{2(-r)^n}{\alpha n^2+(1-\alpha)n}\bigg|=r+\sum_{n=2}^{\infty}\frac{2(-1)^{n-1}r^n}{\alpha n^2+(1-\alpha)n}.
\end{align}
Therefore, in view of \eqref{eee-3.2} and \eqref{ee-33.88}, the distance $ d $ between $ f_{\alpha}(0) $ and the boundary of $ f_{\alpha}(\mathbb{D}) $ is
\begin{equation}\label{e-3.3}
	d(f_{\alpha}(0), \partial f_{\alpha}(\mathbb{D}))=\liminf_{r\rightarrow 1^{-}}	|f_{\alpha}(-r)-f_{\alpha}(0)|=1+\sum_{n=2}^{\infty}\frac{2(-1)^{n-1}}{\alpha n^2+(1-\alpha)n}.
\end{equation}
 A simple computation using \eqref{e-3.22}, \eqref{e-3.23} and \eqref{e-3.3} for $ f=f_{\alpha} $ and $ r>r_k(\alpha) $ shows that 
\begin{align*}
	|z|&+\sum_{n=2}^{\infty}\left(|a_n|+|b_n|\right)|z|^n+P\left(\frac{S_{r}}{\pi}\right)\\&>r_k(\alpha)+\sum_{n=2}^{\infty}\left(|a_n|+|b_n|\right)r^n_k(\alpha)+P\left(\frac{S_{r_k(\alpha)}}{\pi}\right)\\&=r_f(\alpha)+\sum_{n=2}^{\infty}\frac{2r^n_k(\alpha)}{\alpha n^2+(1-\alpha)n}+P\left(r^2_k(\alpha)+\sum_{n=2}^{\infty}\frac{4nr^{2n}_k(\alpha)}{\left(\alpha n^2+(1-\alpha)n\right)^2}\right)\\&= 1+\sum_{n=2}^{\infty}\frac{2(-1)^{n-1}}{\alpha n^2+(1-\alpha)n}\\&=d(f_{\alpha}(0), \partial f_{\alpha}(\mathbb{D})), 
\end{align*}
which shows that $ r_k(\alpha) $ is the best possible. This completes the proof of (i).\vspace{1.5mm}

\noindent (ii) In view of Lemma \ref{lem-1.3} and Lemma \ref{lem-1.5} and \eqref{e-3.21} for $ |z|=r $, we obtain
\begin{align}\label{e-3.24}
	&|f(z)|^m+\sum_{n=2}^{\infty}\left(|a_n|+|b_n|\right)|z|^n+P\left(\frac{S_r}{\pi}\right)\\&\leq \nonumber \left(r+\sum_{n=2}^{\infty}\frac{2r^n}{\alpha n^2+(1-\alpha)n}\right)^m+\sum_{n=2}^{\infty}\frac{2r^n}{\alpha n^2+(1-\alpha)n}\\&\nonumber\quad\quad+P\left(r^2+\sum_{n=2}^{\infty}\frac{4nr^{2n}}{\left(\alpha n^2+(1-\alpha)n\right)^2}\right)^.
\end{align}
A simple computation shows that
\begin{align*}
	&\left(r+\sum_{n=2}^{\infty}\frac{2r^n}{\alpha n^2+(1-\alpha)n}\right)^m+\sum_{n=2}^{\infty}\frac{2r^n}{\alpha n^2+(1-\alpha)n}\\&\nonumber\quad\quad+P\left(r^2+\sum_{n=2}^{\infty}\frac{4nr^{2n}}{\left(\alpha n^2+(1-\alpha)n\right)^2}\right)\leq 1+\sum_{n=2}^{\infty}\frac{2(-1)^{n-1}}{\alpha n^2+(1-\alpha)n}.
\end{align*}
for $ r\leq  r^*_k(m,\alpha)$, where $ r^*_k(m,\alpha) $ is a root of $ J_2(r)=0 $ in $ (0,1) $.\vspace{1.2mm}
Using the similar argument as used in the proof of (i), it can be shown that $ r^*_k(m, \alpha) $ is the unique root of $ J_2(r) =0$. Thus, we have
\begin{align}\label{e-3.25}
	&\left(r^*_k(m,\alpha)+\sum_{n=2}^{\infty}\frac{2\left(r^*_k(m, \alpha)\right)^n}{\alpha n^2+(1-\alpha)n}\right)^m+\sum_{n=2}^{\infty}\frac{2\left(r^*_k(m, \alpha)\right)^n}{\alpha n^2+(1-\alpha)n}\\&\nonumber\quad\quad+P\left(\left(r^*_k(m, \alpha)\right)^2+\sum_{n=2}^{\infty}\frac{4n\left(r^*_k(m, \alpha)\right)^{2n}}{\left(\alpha n^2+(1-\alpha)n\right)^2}\right)= 1+\sum_{n=2}^{\infty}\frac{2(-1)^{n-1}}{\alpha n^2+(1-\alpha)n}.
\end{align}
In view of \eqref{e-3.25} and the arguments used in the part (i) of the theorem, it can be easily shown that $ r^*_k(m, \alpha) $ is the best possible. This completes the proof of (ii).
\end{proof}	
\begin{proof}[\bf Proof of Corollary \ref{cor-2.14}]
	(i) Let $f\in \mathcal{W}^0_{\mathcal{H}}(\alpha) $, then for $ |z|=r $ and $ \alpha=1/2 $, using Lemma \ref{lem-1.3}, we obtain
	\begin{align}\label{e-33.26}
		|z|+\sum_{n=2}^{\infty}\left(|a_n|+|b_n|\right)|z|^n+\frac{S_r}{\pi}\leq 	r^2+r+\sum_{n=2}^{\infty}\frac{4r^n}{n^2+n}+\sum_{n=2}^{\infty}\frac{16nr^{2n}}{(n^2+n)^2}.
	\end{align}
 A simple computation shows that
 \[
 \begin{cases}
 	 \displaystyle\sum_{n=2}^{\infty}\frac{4r^n}{n^2+n}=4-2r+\frac{4}{r}(1-r)\log(1-r),\vspace{1.5mm}\\
 	  \displaystyle\sum_{n=2}^{\infty}\frac{4(-1)^{n-1}}{n^2+n}=2(-3+4\log 2).
 \end{cases}
 \]
On the other hand, we obtain
\begin{align*}
	\sum_{n=2}^{\infty}\frac{16nr^{2n}}{(n^2+n)^2}&=-\sum_{n=2}^{\infty}\frac{16r^{2n+2}}{r^2(n+1)}-\sum_{n=2}^{\infty}\frac{16r^{2n}}{(n+1)^2}+\sum_{n=2}^{\infty}\frac{16r^{2n}}{n}\\&=\frac{16}{r^2}\left(\frac{r^4}{2}+r^2+\log(1-r^2)\right)-\frac{16}{r^2}\left(\frac{r^4}{4}-r^2+{\rm Li}_2(r^2)\right)\\&\quad\quad+16\left(-\log(1-r^2)-r^2\right)\\&=-4r^2+\frac{16}{r^2}\left((1-r^2)\log(1-r^2)\right)-\frac{16}{r^2}{\rm Li}_2(r^2)+32.
\end{align*}
Therefore, we have
\begin{align}
	\label{e-33.27}
	&r^2+r+\sum_{n=2}^{\infty}\frac{4r^n}{n^2+n}+\sum_{n=2}^{\infty}\frac{16nr^{2n}}{(n^2+n)^2}\\&=\nonumber	r^2+r+4-2r+\frac{4}{r}(1-r)\log(1-r)-4r^2+\frac{16(1-r^2)\log(1-r^2)}{r^2}\\&\nonumber\quad\quad-\frac{16{\rm Li}_2(r^2)}{r^2}+32.
\end{align}
Hence, 
\begin{align*}
&r^2+r+4-2r+\frac{4}{r}(1-r)\log(1-r)-4r^2+\frac{16(1-r^2)\log(1-r^2)}{r^2}\\&\nonumber\quad\quad-\frac{16{\rm Li}_2(r^2)}{r^2}+32\leq 2(-3+4\log 2)
\end{align*}
for $ r\leq r_f(1/2) $, where $ r_f(1/2) $ is root of $ F(r)=0 $ in $ (0,1) $.\vspace{1.2mm}

Using the standard argument, we can show that $ F(r) $ has the unique root in $ (0,1) $. Let the root be denoted by $ r_f(1/2) $. Furthermore, a simple computation shows that $ r_f(1/2)\approx 0.600881 $. Hence, we have
\begin{align*}
	&\frac{4}{r_f(1/2)}(1-r_f(1/2))\log (1-r_f(1/2))+\frac{16}{r^2_f(1/2)}(1-r^2_f(1/2))\log(1-r^2_f(1/2))\\&\quad\quad-\frac{16}{r^2_f(1/2)}{\rm Li}_2(r^2_f(1/2))-3r^2_f(1/2)-r_f(1/2)+29+8\log 2=0,
\end{align*}
which is equivalent to
	\begin{align}\label{e-33.28}
	&r^2_f(1/2)+r_f(1/2)+4-2r_f(1/2)+\frac{4}{r_f(1/2)}(1-r_f(1/2))\log(1-r_f(1/2))-4r^2_f(1/2)\\&\nonumber\quad\quad+\frac{16(1-r^2_f(1/2))\log(1-r^2_f(1/2))}{r^2_f(1/2)}-\frac{16{\rm Li}_2(r^2_f(1/2))}{r^2_f(1/2)}+32\\&\nonumber= 1+2(-3+4\log 2).
\end{align}
To show that $ r_f(1/2)\approx 0.600881  $ is the best possible, we consider the function $ f=f_{1/2} $ defined by
\begin{equation}\label{ee-33.20}
	f_{1/2}(z)=z+\sum_{n=2}^{\infty}\frac{4z^n}{n^2+n}.
\end{equation}
In view of \eqref{e-3.3}, an easy computation shows that
\begin{equation}\label{ee-33.21}
	d(f_{1/2}(0), \partial f_{1/2}(\mathbb{D}))=1+\sum_{n=2}^{\infty}\frac{4(-1)^{n-1}}{n^2+n}=1+2(-3+4\log 2).
\end{equation}
Thus, using \eqref{e-33.26}, \eqref{e-33.28} and \eqref{ee-33.21} for $ f=f_{1/2} $ and $ r>r_f(1/2) $, we obtain 
\begin{align*}
	|z|&+\sum_{n=2}^{\infty}\left(|a_n|+|b_n|\right)|z|^n+\frac{S_r}{\pi}\\&>r_f(1/2)+\sum_{n=2}^{\infty}\left(|a_n|+|b_n|\right)\left(r_f(1/2)\right)^n+\frac{S_{r_f(1/2)}}{\pi}\\&=r^2_f(1/2)+r_f(1/2)+4-2r_f(1/2)+\frac{4}{r_f(1/2)}(1-r_f(1/2))\log(1-r_f(1/2))\\&\nonumber\quad\quad-4r^2_f(1/2)+\frac{16(1-r^2_f(1/2))\log(1-r^2_f(1/2))}{r^2_f(1/2)}-\frac{16{\rm Li}_2(r^2_f(1/2))}{r^2_f(1/2)}+32\\&\nonumber= 1+2(-3+4\log 2)\\&=d(f_{1/2}(0), \partial f_{1/2}(\mathbb{D})),
\end{align*}
which shows that $ r_f(1/2) $ is the best possible. This completes the proof of (i).\vspace{1.5mm}

\noindent (ii) Using Lemma \ref{lem-1.3}, for $ |z|=r $, a simple computation shows that
\begin{align}\label{e-33.29}
	&|f(z)|+\sum_{n=2}^{\infty}\left(|a_n|+|b_n|\right)|z|^n+\frac{S_r}{\pi}\\&\nonumber\leq r+\sum_{n=2}^{\infty}\frac{4r^n}{n^2+n}+\left(r^2+\sum_{n=2}^{\infty}\frac{16nr^{2n}}{n^2+n}\right)\\&\nonumber=-3r^2-3r+\frac{8}{r}(1-r)\log(1-r)+\frac{16(1-r^2)\log(1-r^2)}{r^2}-\frac{16{\rm Li}_2(r^2)}{r^2}+40.
\end{align}
It is easy to see that
\begin{align*}
	&-3r^2-3r+\frac{8}{r}(1-r)\log(1-r)+\frac{16(1-r^2)\log(1-r^2)}{r^2}-\frac{16{\rm Li}_2(r^2)}{r^2}+40\\&\leq 1+2(-3+4\log 2)
\end{align*}
for $ r\leq r^*_f(1/2) $, where $ r^*_f(1/2) $ is a root of $ T(r)=0 $ in $ (0,1) $.\vspace{1.2mm}

By a similar argument as used in part (i), it can be shown that $ T(r) $ has the unique root in $ (0,1) $. Let $ r^*_f(1/2) $ be the root of $ T(r) $ in $ (0, 1) $. Then, we have $ r^*_f(1/2)\approx 0.302059 $. Therefore, we must have
\begin{align*}
	&-3\left(r^*_f(1/2)\right)^2-3r^*_f(1/2)+\frac{8}{r^*_f(1/2)}(1-r^*_f(1/2))\log(1-r^*_f(1/2))\\&\quad\quad+\frac{16(1-\left(r^*_f(1/2)\right)^2)\log(1-\left(r^*_f(1/2)\right)^2)}{\left(r^*_f(1/2)\right)^2}-\frac{16{\rm Li}_2(\left(r^*_f(1/2)\right)^2)}{\left(r^*_f(1/2)\right)^2}\\&\quad\quad+45-8\log 2=0
\end{align*}
and this is equivalent to
\begin{align}\label{e-33.30}
	&-3\left(r^*_f(1/2)\right)^2-3r^*_f(1/2)+\frac{8}{r^*_f(1/2)}(1-r^*_f(1/2))\log(1-r^*_f(1/2))\\&\nonumber\quad\quad+\frac{16(1-\left(r^*_f(1/2)\right)^2)\log(1-\left(r^*_f(1/2)\right)^2)}{\left(r^*_f(1/2)\right)^2}-\frac{16{\rm Li}_2(\left(r^*_f(1/2)\right)^2)}{\left(r^*_f(1/2)\right)^2}+40\\&\nonumber=1+2(-3+4\log 2).
\end{align}
In order to show that $ r^*_f(1/2)\approx 0.302059 $ is the best possible, we consider the function $ f=f_{1/2} $ defined by \eqref{ee-33.20}. Then in view of \eqref{ee-33.21}, \eqref{e-33.29} and \eqref{e-33.30}, for $ f=f_{1/2} $ and $ r>r^*_f(1/2) $, by a similar arguments a used in part (i) shows that
\begin{align*}
|f(z)|+\sum_{n=2}^{\infty}\left(|a_n|+|b_n|\right)|z|^n+\frac{S_{r}}{\pi}> d(f_{1/2}(0), \partial f_{1/2}(\mathbb{D})).
\end{align*}
This shows that $ r^*_f(1/2) $ is the best possible. This completes the proof of (ii).
\end{proof}	
\begin{proof}[\bf Proof of Theorem \ref{th-3.8}]
	(i) In view part-(i) of Lemma \ref{lem-3.2} and \eqref{e-3.18}, we see that $ S_r/\pi\leq r^2/(1-r^2)^2 $. Also, $ M_f(r)\leq\sum_{n=1}^{\infty}r^n=r/(1-r) $. Therefore, a simple computation shows that
	\begin{align}\label{e-6.21}
		M_f(r)+P\left(\frac{S_r}{\pi}\right)\leq \frac{r}{1-r}+\frac{\lambda_1 r^2}{(1-r^2)^2}+\cdots+\frac{\lambda_k r^{2k}}{(1-r^2)^{2k}}\leq \frac{1}{2}
	\end{align}
for $ |z|=r\leq r_k $, where $ r_k $ is a root in $ (0, 1) $ of the equation $ F(r)=0 $, where
\begin{align*}
	F(r):=\frac{r}{1-r}+\frac{\lambda_1 r^2}{(1-r^2)^2}+\cdots+\frac{\lambda_k r^{2k}}{(1-r^2)^{2k}}-\frac{1}{2}.
\end{align*}
The existence of $ r_k $ in $ (0, 1) $ is confirmed by the fact that $ F $ is a real valued differential function in $ (0, 1) $ with the properties $ F(0)=-\frac{1}{2}<0 $ and $ \lim\limits_{r\rightarrow 1^{-}} F(r)=+\infty $. On other hand, a simple computation shows that
\begin{align*}
	\frac{d}{dr}(F(r))=\frac{1}{(1-r)^2}+2\lambda_1 g(r)g^{\prime}(r)+\cdots+2k\lambda_kg^{2k-1}(r)g^{\prime}(r)>0\; \mbox{for}\; r\in (0, 1),
\end{align*}
where
\begin{align*}
	g(r):=\frac{r}{1-r^2}\;\mbox{and hence}\; g^{\prime}(r)=\frac{1+r^2}{(1-r^2)^2}\;\; r\in (0, 1).
\end{align*}
Thus, $ F $ being monotone increasing, the root $ r_k $ is unique. Thus, we have
\begin{align}\label{e-6.22}
	\frac{r_k}{1-r_k}+\frac{\lambda_1 r_k^2}{(1-r_k^2)^2}+\cdots+\frac{\lambda_k r_k^{2k}}{(1-r_k^2)^{2k}}=\frac{1}{2}.
\end{align}
In view of part-(iii) of Lemma \ref{lem-3.2}, we see that
\begin{align}\label{e-6.23}
	d(f(0),\partial f(\mathbb{D}))=\liminf_{|z|=r\rightarrow 1^{-}} |f(z)-f(0)|\geq \liminf_{|z|=r\rightarrow 1^{-}}\frac{|z|}{1+|z|}=\frac{1}{2}.
\end{align}
Therefore, it follows from \eqref{e-6.21}, \eqref{e-6.22}and \eqref{e-6.23} that
\begin{align*}
		M_f(r)+P\left(\frac{S_r}{\pi}\right)\leq d(f(0),\partial f(\mathbb{D}))
\end{align*}
for $ |z|=r\leq r_k $.\vspace{1.2mm}

In order to show that $ r_k $ is best possible, we consider a suitable rotation of the analytic mapping $ f^{*}(z)=z/(1-z) $ in $ \mathbb{D} $. It can be shown that $ M_{f^{*}}(r)={r}/{(1-r)} $, $ S_r/\pi=r^2/(1-r^2)^2 $ and $ d(f^{*}(0),\partial f^{*}(\mathbb{D}))=1/2 $. For this function $ f^{*} $ and $ r>r_k $, in view of \eqref{e-6.22}, we see that
\begin{align*}
	M_{f^{*}}(r)+P\left(\frac{S_r}{\pi}\right)&=\frac{r}{1-r}+\frac{\lambda_1 r^2}{(1-r^2)^2}+\cdots+\frac{\lambda_k r^{2k}}{(1-r^2)^{2k}}\\&>\frac{r_k}{1-r_k}+\frac{\lambda_1 r_k^2}{(1-r_k^2)^2}+\cdots+\frac{\lambda_k r_k^{2k}}{(1-r_k^2)^{2k}}&=\frac{1}{2}&=d(f^{*}(0),\partial f^{*}(\mathbb{D})),
\end{align*}
which shows that $ r_k $ is best possible.\vspace{1.2mm}

	(ii) By part-(ii) of Lemma \ref{lem-3.2} and \eqref{e-3.18}, we see that 
	\begin{align*}
		\frac{S_r}{\pi}\leq\sum_{n=1}^{\infty} n^3r^{2n}=\frac{r^6+4r^4+r^2}{(r^2-1)^2}.
	\end{align*} 
Note that $ M_f(r)\leq r/(1-r)^2 $. Therefore, a simple computation shows that
\begin{align}\label{e-6.24}
	M_f(r)+P\left(\frac{S_r}{\pi}\right)\leq \frac{r}{(1-r)^2}+\frac{\lambda_1\left(r^6+4r^4+r^2\right)}{(r^2-1)^2}+\cdots+\frac{\lambda_k \left(r^6+4r^4+r^2\right)^{k}}{(r^2-1)^{2k}}\leq \frac{1}{4}
\end{align}
for $ |z|=r\leq r^{*}_k $, where $ r_k $ is a root in $ (0, 1) $ of the equation $ G(r)=0 $, where
\begin{align*}
	G(r):=\frac{r}{(1-r)^2}+\frac{\lambda_1\left(r^6+4r^4+r^2\right)}{(r^2-1)^2}+\cdots+\frac{\lambda_k \left(r^6+4r^4+r^2\right)^{k}}{(r^2-1)^{2k}}-\frac{1}{4}.
\end{align*}
By the similar argument used in part-(i), it can be shown that the root $ r^{*}_k $ in $ (0, 1) $ is unique. In view of part-(iv) of Lemma \ref{lem-3.2}, we now see that
\begin{align}\label{e-6.25}
	d(f(0),\partial f(\mathbb{D}))=\liminf_{|z|=r\rightarrow 1^{-}} |f(z)-f(0)|\geq \liminf_{|z|=r\rightarrow 1^{-}}\frac{|z|}{\left(1+|z|\right)^2}=\frac{1}{4}.
\end{align}
Therefore, it follows from \eqref{e-6.24} and \eqref{e-6.25} that
\begin{align*}
	M_f(r)+P\left(\frac{S_r}{\pi}\right)\leq d(f(0),\partial f(\mathbb{D}))
\end{align*}
for $ |z|=r\leq r^{*}_k $.\vspace{1.2mm}

To show that $ r^{*}_k $ is best possible, we consider a suitable rotation of the analytic mapping $ f^{**}(z)=z/(1-z)^2 $ in $ \mathbb{D} $. It can be shown that $ M_{f^{**}}(r)={r}/{(1-r)^2} $, $ S_r/\pi={\left(r^6+4r^4+r^2\right)}/{(r^2-1)^2} $ and $ d(f^{**}(0),\partial f^{**}(\mathbb{D}))=1/4 $. For this function $ f^{**} $ and $ r>r^{*}_k $, by the similar argument used in part-(i) of the theorem, it can be easily shown that $ r^{*}_k $ is best possible.
\end{proof}
\begin{proof}[\bf Proof of Theorem \ref{th-3.9}]
	(i) In view part-(i) of Lemma \ref{lem-3.2}, \eqref{e-33.8} and \eqref{e-3.18}, we see that
	\begin{align*}
		&M_f(r)+P\left(\frac{S_r}{\pi-S_r}\right)\\&\nonumber\leq \frac{r}{1-r}+\dfrac{\lambda_k r^{2k}}{((1-r^2)^2-r^2)^k}+\cdots+\dfrac{\lambda_1r^2}{(1-r^2)^2-r^2}\\&\nonumber\leq \frac{1}{2}
	\end{align*}
for $ |z|=r\leq R_k $, where $ R_k $ is a root of $ F_1(r)=0 $ in $ (0, 1) $, where
\begin{align*}
	F_1(r)=\frac{r}{1-r}+\dfrac{\lambda_k r^{2k}}{((1-r^2)^2-r^2)^k}+\cdots+\dfrac{\lambda_1r^2}{(1-r^2)^2-r^2}-\frac{1}{2}.
\end{align*}
By the similar argument used in the part (i) of Theorem \ref{th-3.8}, it can be shown that 
\begin{align*}
	M_f(r)+P\left(\frac{S_r}{\pi-S_r}\right)\leq d(f(0), \partial f(\mathbb{D}))
\end{align*}
for $ |z|=r\leq R_k $ where $ R_k $ is the smallest root in $ (0, 1) $ of the equation 
\begin{align*}
	\frac{r}{1-r}+\dfrac{\lambda_k r^{2k}}{((1-r^2)^2-r^2)^k}+\cdots+\dfrac{\lambda_1r^2}{(1-r^2)^2-r^2}=\frac{1}{2}
\end{align*}
and $ R_k $ is best possible. Hence, we omit the details.\vspace{1.2mm}

(ii) (i) In view part-(ii) of Lemma \ref{lem-3.2}, \eqref{e-33.9}  and \eqref{e-3.18}, a simple computation shows that
\begin{align*}
	&M_f(r)+P\left(\frac{S_r}{\pi-S_r}\right)\\&\leq \nonumber	\frac{r}{(1-r)^2}+\lambda_k\left(\dfrac{r^6+4r^4+r^2}{(r^2-1)^4-(r^6+4r^4+r^2)}\right)^k\\&\nonumber\quad+\cdots+\lambda_1\left(\dfrac{r^6+4r^4+r^2}{(r^2-1)^4-(r^6+4r^4+r^2)}\right)\\&\nonumber\leq \frac{1}{4}
\end{align*}
for $ |z|=r\leq R^{*}_k $, where $ R^{*}_k $ is a root of $ F_2(r)=0 $ in $ (0, 1) $, where
\begin{align*}
	F_1(r)=	&\frac{r}{(1-r)^2}+\lambda_k\left(\dfrac{r^6+4r^4+r^2}{(r^2-1)^4-(r^6+4r^4+r^2)}\right)^k\\&+\cdots+\lambda_1\left(\dfrac{r^6+4r^4+r^2}{(r^2-1)^4-(r^6+4r^4+r^2)}\right)-\frac{1}{4}.
\end{align*}
By the similar argument used in the part (ii) of Theorem \ref{th-3.8}, it can be easily shown that 
\begin{align*}
	M_f(r)+P\left(\frac{S_r}{\pi-S_r}\right)\leq d(f(0), \partial f(\mathbb{D}))
\end{align*}
for $ |z|=r\leq R^{*}_k $ where $ R^{*}_k $ is the smallest root in $ (0, 1) $ of the equation 
\begin{align*}
		&\frac{r}{(1-r)^2}+\lambda_k\left(\dfrac{r^6+4r^4+r^2}{(r^2-1)^4-(r^6+4r^4+r^2)}\right)^k\\&+\cdots+\lambda_1\left(\dfrac{r^6+4r^4+r^2}{(r^2-1)^4-(r^6+4r^4+r^2)}\right)=\frac{1}{4}.
\end{align*}
and $ R^{*}_k $ is best possible. This completes the proof.\vspace{1.2mm}
\end{proof}
\vspace{2mm}
\noindent\textbf{Acknowledgment:} The second author is supported by SERB-CRG, India.\\

\vspace{2.5mm}

\noindent\textbf{ Compliance of Ethical Standards}\\

\noindent\textbf{ Conflict of interest} The authors declare that there is no conflict  of interest regarding the publication of this paper. \\

\noindent\textbf{ Data availability statement}  Data sharing not applicable to this article as no datasets were generated or analysed during the current study.

\end{document}